\documentclass[a4paper,10pt]{article}
\everymath{\displaystyle}
\usepackage{amssymb,amsmath,amsthm}
\usepackage{graphicx}

\def\N{\mathbb N}
\def\R{\mathbb R}
\def\T{\mathbb T}

\usepackage{amsfonts}
\usepackage{hyperref}
\usepackage{latexsym}
\usepackage{amsthm}
\usepackage{amsmath}
\usepackage{amssymb}
\usepackage{amsopn}
\usepackage{geometry}
\usepackage{amscd}
\usepackage{mathrsfs}
\usepackage{dsfont}
\usepackage[english]{babel}
\usepackage{caption}
\font\teneufm=eufm10
\font\seveneufm=eufm7
\font\fiveeufm=eufm5
\newfam\eufmfam
\textfont\eufmfam=\teneufm
\scriptfont\eufmfam=\seveneufm
\scriptscriptfont\eufmfam=\fiveeufm
\def\eufm#1{{\fam\eufmfam\relax#1}}

\def\vgot{{\eufm v}}

\newcommand{\Om}{\Omega}

\newcommand{\aaa}{{\mathtt{a}}}

\newtheorem{theorem}{Theorem}[section]
\newtheorem{proposition}{Proposition}[section]

\newtheorem{remark}{Remark}[section]
\newtheorem{remarks}{Remark}[section]
\newtheorem{definition}{Definition}[section]
\newcommand{\be}{\begin{equation}}
\newcommand{\ee}{\end{equation}}

\newcommand{\om}{\omega}
\newcommand{\e}{\varepsilon}

\renewcommand{\a }{\alpha }
\renewcommand{\b }{\beta }
\newcommand{\s }{\sigma }
\newcommand{\ii }{{\rm i} }

\newcommand{\g }{\gamma}

\renewcommand{\l }{\lambda }

\renewcommand{\o }{\omega }
\renewcommand{\O }{\Omega }
\newcommand{\C}{\mathbb{C}}
\newcommand{\Z}{\mathbb{Z}}

\newcommand{\mm}{{\rm m}}

\newcommand{\norma}{\|}

\newcommand{\OO}{\O}
\newcommand{\oo}{\mathtt{\o}}

\usepackage{color}

\title{\bf Existence and stability of quasi-periodic solutions for
derivative wave equations}

\author{Massimiliano Berti, Luca Biasco, Michela Procesi}

\date{}

\begin{document}

\maketitle

\noindent
{\bf Abstract:} In this note we present the new KAM result in \cite{BBP3}  which 
proves the existence of Cantor families of small amplitude, analytic,
quasi-periodic solutions  of derivative wave equations,
with zero Lyapunov exponents and whose linearized equation is reducible to constant coefficients.
In turn, this result is derived by an abstract KAM theorem for   infinite dimensional  reversible dynamical systems.
\\[1.5mm]
{\sc 2000AMS subject classification}: 37K55, 35L05.


\section{Introduction}

In the last years many progresses have been obtained concerning  KAM theory
for nonlinear PDEs, 
since  the pioneering works of  Kuksin
\cite{Ku} and Wayne \cite{W1} for  $1$-$d$ semilinear wave (NLW) and
Schr\"odinger (NLS) equations under Dirichlet boundary conditions, 
see  \cite{KP} and \cite{Po3} for further developments.
The approach of these papers consists  in generating iteratively
a sequence of symplectic changes of variables 
which bring the Hamiltonian into a  constant coefficients ($=$reducible) normal form  with an elliptic  ($=$linearly stable)
invariant torus at the origin.
Such a  torus is filled by quasi-periodic solutions with zero Lyapunov exponents.
This procedure requires
 to solve, at each step, constant-coefficients linear ``homological equations"
by imposing the ``second order Melnikov" non-resonance conditions.
Unfortunately these (infinitely many)
conditions are violated already for periodic boundary conditions. 

In this case, existence of quasi-periodic solutions for semilinear $ 1d $-NLW and NLS equations,
was first proved by Bourgain \cite{B1} by extending the Lyapunov-Schmidt decomposition and the 
Newton approach introduced by Craig-Wayne \cite{CW} for periodic solutions.
Its main advantage is to require only the  ``first order Melnikov"
non-resonance conditions (the minimal assumptions)  for solving the
homological equations. 
It has allowed 
Bourgain  to prove  \cite{B3}, \cite{B5}  also
the existence of
quasi-periodic solutions  for NLW and NLS (with Fourier multipliers) in any space dimension,
see also the recent extensions  in Berti-Bolle \cite{BBo}, \cite{BB12}. 
The main drawback of this approach is that the
homological equations are  linear PDEs with non-constant coefficients.
Translated in the KAM language this implies
a non-reducible normal form around the torus and then a lack of  information about the
stability of the quasi-periodic solutions.
Later on, existence of reducible elliptic  tori was proved 
by   Eliasson-Kuksin \cite{EK} for NLS  (with Fourier multipliers) in any space dimension,
see also Procesi-Xu \cite{PX}. 

\smallskip

 A challenging frontier concerns PDEs with unbounded 
 nonlinearities, i.e. containing derivatives. 
In this direction KAM theory has been extended for perturbed KdV  equations
 by Kuksin \cite{K2},
Kappeler-P\"oschel \cite{KaP},  and,
for the $1d$-derivative NLS (DNLS) and Benjiamin-Ono equations,  by Liu-Yuan \cite{LY}.
We remark that the KAM proof   is more delicate for DNLS and Benjiamin-Ono, because these equations
 are less ``dispersive" than KdV, i.e.
  the eigenvalues of the principal part of the differential operator   grow only quadratically  at infinity, and not cubically as  for KdV.
This difficulty is reflected in the fact that
the quasi-periodic solutions  in \cite{K2}, \cite{KaP} are analytic, while those of    \cite{LY} are only $ C^\infty $.
 Actually, for the applicability of these KAM schemes,
the more dispersive the equation is,
the more derivatives in the nonlinearity  can be supported. 
The limit case 
of the derivative nonlinear wave equation (DNLW) -which is not dispersive at all-  is excluded  by these approaches.

\smallskip

In this note we present the 
KAM theory developed in \cite{BBP3} which proves 
existence and stability of small amplitude analytic quasi-periodic solutions of the derivative wave equations. Such PDEs
are not  Hamiltonian, but may have a reversible structure, that we shall exploit. 

\smallskip

All the previous results 
concern  {\it Hamiltonian} PDEs.
It was however remarked by Bourgain that the construction of periodic and
quasi-periodic solutions,  using 
the  Newton iteration method  of Craig-Wayne \cite{CW}, 
is a-priori not restricted to Hamiltonian systems. This  approach
appears  as a general implicit function type result, in large part
independent of the Hamiltonian character of the equations. For
example in  \cite{B2} Bourgain  proved the existence of periodic
solutions for the non-Hamiltonian derivative wave equation
\be\label{yttyxx}
{\mathtt y}_{tt} - {\mathtt y}_{{\mathtt
x}{\mathtt x}} + \mm {\mathtt y} + {\mathtt y}_t^2  = 0 \, , \quad
\mm \neq 0 \, , \quad {\mathtt x} \in \T := \R \slash (2 \pi \Z )\, .
\ee
Actually also KAM theory is not only  Hamiltonian in nature, but may be formulated 
for general vector fields, as
realized in the seminal work of Moser \cite{Mo67}.
This paper, in particular,  
started the analysis of reversible KAM theory for finite
dimensional systems, later extended by Arnold \cite{Arn} and
Servyuk \cite{Se1}. 
The reversibility property implies that 
the averages over the fast angles of some components of the vector field are zero, thus  
removing the ``secular drifts" of the actions which are  
incompatible  with a quasi-periodic behavior of the solutions.

\smallskip

Recently, Zhang-Gao-Yuan \cite{ZGY}  have proved the
existence of  $ C^\infty$-quasi periodic solutions for the derivative NLS
equation 
$$
\ii u_t + u_{{\mathtt x}{\mathtt x}} +
|u_{\mathtt x}|^2 u = 0 
$$ 
with Dirichlet boundary conditions.
Such equation  is reversible, but  not Hamiltonian. 
The  result  \cite{ZGY}  is proved adapting the KAM scheme
developed for the Hamiltonian DNLS in Liu-Yuan \cite{LY}. 
The derivative nonlinear wave equation (DNLW), which is
not dispersive, is excluded. 

\smallskip

In the recent paper \cite{BBP1} we have extended KAM theory to deal
with Hamiltonian  derivative wave equations like
$$
{\mathtt y}_{tt} - {\mathtt y}_{{\mathtt x}{\mathtt x}} + \mm
{\mathtt y} + f(D {\mathtt y})=0 \, , \quad \mm >  0 \, , \quad D
:=\sqrt{-\partial_{{\mathtt x}{\mathtt x}}+\mm} \, , \quad
{\mathtt x} \in \T \, .
$$
This kind of  Hamiltonian pseudo-differential equations has been
introduced by Bourgain  \cite{B1} and Craig  \cite{C} as
models to study the effect of derivatives versus dispersive
phenomena. The key of  \cite{BBP1} is the proof of 
the  first order asymptotic expansion of the perturbed normal
frequencies, obtained using the notion of quasi-T\"oplitz function. This concept was introduced by Procesi-Xu 
\cite{PX} and  is connected to the T\"oplitz-Lipschitz property in Eliasson-Kuksin \cite{EK}. 
Of course we could
not deal in \cite{BBP1}  with the derivative wave equation, which is not 
Hamiltonian.

\smallskip

In \cite{BBP3} we develop KAM theory for a class of reversible derivative wave equations 
\be\label{DNLW} 
{\mathtt
y}_{tt} - {\mathtt y}_{{\mathtt x}{\mathtt x}} + \mm {\mathtt y} =
g({\mathtt x}, {\mathtt y}, {\mathtt y}_{\mathtt x}, {\mathtt y}_t) \, , \quad {\mathtt x} \in \T  \, ,
\ee
implying the {\it existence} and the {\it stability} of {\it analytic} quasi-periodic solutions,
see Theorem \ref{thm:DNLW}.
Note that the nonlinearity in \eqref{DNLW} has an explicit $ \mathtt x $-dependence (unlike \cite{BBP1}).
 The search for periodic/quasi--periodic solutions for derivative wave equations
 is a natural question,  which was pointed out, for instance, 
  by Craig  \cite{C} as an important open problem (see section 7.3 of \cite{C}).
\smallskip

Clearly we can not  expect  an existence result for any nonlinearity. 
For example,  \eqref{DNLW} with the nonlinear friction term $ g = {\mathtt y}_t^3 $ has no smooth
periodic/quasi-periodic solutions except the constants, see Proposition \ref{nonext}.
This case is ruled out by assuming the  condition
\be\label{geven}
g({\mathtt x}, {\mathtt y}, {\mathtt
y}_{\mathtt x}, - {\mathtt v} ) = g({\mathtt x}, {\mathtt y},
{\mathtt y}_{\mathtt x},  {\mathtt v} ) \, , 
\ee
satisfied, for example, by  \eqref{yttyxx}.
Under condition \eqref{geven} the equation \eqref{DNLW} 
is {\it reversible}, namely the associated  first order system
\be\label{firstor}
{\mathtt y}_t =  {\mathtt v}  \, , \quad
{\mathtt v} _{t} = {\mathtt y}_{{\mathtt x}{\mathtt x}} - \mm
{\mathtt y} + g({\mathtt x}, {\mathtt y}, {\mathtt y}_{\mathtt x},
{\mathtt v}  ) \, ,
\ee
is reversible with respect to the
involution
\be\label{invol1} 
S ( {\mathtt y},  {\mathtt v}  ) := (
{\mathtt y}, -  {\mathtt v}  ) \, \, , \quad S^2 = I \, .
\ee
 Reversibility is an important property in order to allow the existence of
periodic/quasi-periodic solutions, albeit
not sufficient.  For example, the reversible equation
$$ 
{\mathtt y}_{tt} - {\mathtt y}_{{\mathtt x}{\mathtt x}}   
= {\mathtt y}_{\mathtt x}^3 \, ,
\quad {\mathtt x} \in {\T } \, ,
$$ (proposed in \cite{C}, page 89)
has no smooth periodic/quasi-periodic solutions except the constants, see Proposition \ref{nonex}.
In order to find quasi-periodic solutions we also require the
parity assumption
\be\label{REALITY}
g( - {\mathtt x}, {\mathtt y}, - {\mathtt y}_{\mathtt x} ,  {\mathtt v}  ) =
 g( {\mathtt x}, {\mathtt y},  {\mathtt y}_{\mathtt x} ,  {\mathtt v}  ) \, ,
\ee
which  rules out  nonlinearities like $ {\mathtt y}_{\mathtt x}^3 $.
Actually, for the wave equation  \eqref{DNLW} the role of
the time and space variables $ (t, {\mathtt x}) $ is highly
symmetric. Then, considering $ {\mathtt x} $ ``as time" (spatial
dynamics idea) the term $ {\mathtt y}_{\mathtt x}^3 $ is a
friction and condition \eqref{REALITY} is the corresponding
reversibility condition. 

After  Theorem \ref{thm:DNLW} we shall further comment on the assumptions.

\smallskip

Before stating our main results, we mention the classical
bifurcation theorems of Rabinowitz \cite{R68} about periodic solutions  (with rational periods)
of dissipative forced derivative wave equations
$$
{\mathtt y}_{tt} - {\mathtt y}_{{\mathtt x}{\mathtt x}}  + \a
{\mathtt y}_t + \e F({\mathtt x},t, {\mathtt y}_{\mathtt x},
{\mathtt y}_t ) = 0 \, ,  \quad {\mathtt x} \in [0,\pi]
$$
with Dirichlet boundary conditions,
and 
for fully-non-linear forced wave equations
$$
{\mathtt y}_{tt} - {\mathtt y}_{{\mathtt
x}{\mathtt x}}  + \a {\mathtt y}_t + \e F({\mathtt x},t, {\mathtt
y}_{\mathtt x}, {\mathtt y}_t , {\mathtt y}_{tt}, {\mathtt
y}_{t{\mathtt x}}, {\mathtt y}_{{\mathtt x}{\mathtt x}}) = 0 \, \, , \quad {\mathtt x} \in [0,\pi] \, . 
$$
This latter result is quite subtle because, from the point of view
of the initial value problem, it is uncertain whether a solution
can exist for more than a finite time due to the formation of
shocks. Here  the presence of the dissipation $ \a \neq 0 $
allows the existence of a periodic solution.

\smallskip

Finally, concerning quasi-linear wave equations we mention the Birkhoff normal form results of Delort 
\cite{De1} (and references therein), 
which imply long time existence for solutions with small initial data. To our knowledge,  
these are the only results of this type on compact manifolds.  
For quasi-linear wave equations in $ \R^d $ there is a huge literature, 
since the nonlinear  effects of derivatives may be controlled by dispersion. 

\subsection{Main results} 

We consider derivative wave equations \eqref{DNLW} where $ \mm > 0 $, the nonlinearity
$$
g : \T \times {\cal U} \to \R \, , \quad  {\cal U} \subset \R^3 \ {\rm open \, neighborhood \ of } \ 0  \, ,
$$
is real analytic and satisfies the ``reversibility" and ``parity" assumptions \eqref{geven}, \eqref{REALITY}.
Moreover $ g $ vanishes at least quadratically at $ ({\mathtt y},{\mathtt
y}_{\mathtt x},  \mathtt v ) = (0,0,0 ) $,
namely
\be\label{quadratic}
g({\mathtt x}, 0,0,0)  = (\partial_{\mathtt y} g)({\mathtt x},
0,0,0 ) = (\partial_{{\mathtt y}_{\mathtt x}} g)({\mathtt x},
0,0,0 ) = (\partial_{\mathtt v} g)({\mathtt x}, 0,0,0 ) = 0 \, .
\ee
In addition we assume a ``non-degeneracy" condition on the leading order term of the nonlinearity
(in order to verify the usual ``twist" hypotheses required in KAM theory). 
For definiteness, we have developed all the calculations for  
\be\label{gleading}
g = {\mathtt y} {\mathtt y}_{\mathtt x}^2 + {\rm h.o.t}\, . 
\ee
Because of \eqref{geven}, it is natural to look for ``reversible" quasi-periodic
solutions such that 
$$  
S ({\mathtt y}, {\mathtt v}) (t ) = ({\mathtt y}, {\mathtt v}) (- t ) \, ,  
$$ namely such that $ {\mathtt y}(t, {\mathtt x} ) $ is
even and $  {\mathtt v} (t,{\mathtt x}) $ is odd in time.
Moreover, because of  \eqref{REALITY} the phase space of functions 
even in $ \mathtt x $, 
\be\label{evenx} 
({\mathtt y}, {\mathtt v} )(- {\mathtt x}) =
({\mathtt y}, {\mathtt v} )({\mathtt x}) \, , \  \,  \forall
{\mathtt x \in \T } \, ,
\ee
is invariant under the  flow evolution of \eqref{firstor} and 
it is natural to study  the dynamics on this subspace 
(standing waves).
Note, in particular, that $ {\mathtt y} $ 
satisfies  the {\it Neumann} boundary conditions  $  {\mathtt
y}_{\mathtt x}(t,0) =  {\mathtt y}_{\mathtt x} (t, \pi )  = 0 $.

Summarizing we  look for reversible quasi-periodic {\it standing
wave} solutions  of \eqref{DNLW}, namely satisfying
\be\label{simmet}
{\mathtt y}(t,{\mathtt x}) = {\mathtt
y}(t,-{\mathtt x}) \, , \ \forall t \, , \quad {\mathtt
y}(-t,{\mathtt x}) = {\mathtt y}(t,{\mathtt x}) \, , \  \forall
{\mathtt x \in \T } \, .
\ee
For every choice of the {\sl tangential} sites $ {\cal I}^+ \subset \N
\setminus \{0 \} $,  the linear Klein-Gordon equation
$$ 
{\mathtt y}_{tt} - {\mathtt y}_{{\mathtt x}{\mathtt
x}} + \mm {\mathtt y} = 0  \, , \quad {\mathtt x} \in \T  \, , 
$$
possesses  the family of quasi-periodic standing wave solutions
\be\label{ylineare} 
{\mathtt y} = \sum_{j \in {\cal I}^+}
\frac{\sqrt{8 \xi_j}}{\l_j} \cos ( \l_j  \, t ) \cos (j {\mathtt
x}) \, , \quad \l_j := \sqrt{j^2 + \mm} \, , 
\ee 
parametrized by the amplitudes $ \xi_j \in \R_+ $.

\begin{theorem}\label{thm:DNLW} {\bf (\cite{BBP3})}
For all $ \mm > 0 $, for every choice of finitely many tangential sites
$ {\cal I}^+  \subset \N \setminus \{ 0 \} $, the DNLW equation
\eqref{DNLW} with a real analytic nonlinearity  satisfying \eqref{geven}, \eqref{REALITY}, \eqref{quadratic},
 \eqref{gleading} 
admits small-amplitude, analytic (both in $ t $ and $ \mathtt x $),
quasi-periodic solutions 
\be\label{solutions} 
{\mathtt y}  =
\sum_{j \in {\cal I}^+} \frac{\sqrt{ 8 \xi_j}}{\l_j} \,  \cos (
\om_j^\infty (\xi) \, t ) \cos ( j {\mathtt x} ) + o( \sqrt{\xi}
), \quad
 \om^\infty_j (\xi) \stackrel{\xi \to 0}\approx \sqrt{j^2 + \mm}
 \ee
satisfying \eqref{simmet}, 
for all values of the parameters $\xi$ in a Cantor-like set 
with asymptotic density $ 1$ at $ \xi = 0 $. Moreover the solutions have {\sl zero Lyapunov exponents} and
the linearized equations can be {\sl reduced to constant
coefficients}. The term $ o( \sqrt{\xi} ) $ in \eqref{solutions}
is small in some analytic norm.
\end{theorem}

This theorem  completely
answers the question posed by Craig in \cite{C} 
of developing a general theory for quasi-periodic solutions for
reversible derivative wave equations.
With respect to Bourgain \cite{B2}, we prove the existence of
quasi-periodic solutions (not only periodic) as well as  a stable
KAM normal form nearby.

\smallskip

Let us comment on 
 the hypothesis of Theorem \ref{thm:DNLW}.

\begin{enumerate}
\item {\bf Reversibility and Parity.} As already said,
 the ``reversibility" and  ``parity" assumptions \eqref{geven}, \eqref{REALITY},   rule out
 nonlinearities like $ {\mathtt y}_t^3 $ and $ {\mathtt y}_{\mathtt x}^3 $ for which periodic/quasi-periodic  solutions of \eqref{DNLW}
do not exist. We generalize 
these non-existence results in 
Propositions \ref{nonex}, \ref{nonext}. 
\item
{\bf Mass $ \bf \mm > 0 $.}
The assumption on the mass $ \mm \neq  0 $ is, in general,
necessary.
When $ \mm = 0 $, a well known example of Fritz John (see \eqref{FJ}) 
shows that \eqref{yttyxx} has no smooth solutions for all times except the
constants. In Proposition \ref{FritzJ} we give other
non-existence results of periodic/quasi-periodic solutions for DNLW equations 
satisfying both \eqref{geven}, \eqref{REALITY}, but  with mass $ \mm = 0 $. 
For the KAM construction, the mass $ \mm > 0 $ is used in the Birkhoff normal form step. 
If the mass $ \mm <  0 $ then the Sturm-Liouville operator 
$ - \partial_{{\mathtt x}{\mathtt x}} + \mm $ may possess finitely many
negative eigenvalues and one should expect the existence of partially hyperbolic tori.
\item
{\bf $ \bf \mathtt x $-dependence.} 
The nonlinearity $ g $ in \eqref{DNLW} may explicitly depend on the space
variable $ \mathtt x $, i.e. this equations are
{\sc not} invariant under $ \mathtt x $ translations.
This is an important novelty  with respect to the KAM theorem in
\cite{BBP1} which used the conservation of momentum. 
The key idea
is the introduction of the $ \mathtt a$-weighted majorant norm for
vector fields (see \eqref{aweight}) which  penalizes the
``high-momentum monomials", see \eqref{penalize}.  
\item
{\bf Twist.}  
We have developed all the calculations for
the cubic leading term 
$ g = {\mathtt y} {\mathtt y}_{\mathtt x}^2 + {\rm h.o.t. } $. 
In this case the third order Birkhoff normal form of the PDE \eqref{DNLW} 
turns out to be (partially) integrable and the frequency-to-action map 
is invertible. This is the so called  "twist-condition" in KAM theory.  
It could be interesting to classify the allowed nonlinearities. 
For example, among the cubic nonlinearities, 
we already know that 
for $ {\mathtt y}_{\mathtt x}^3  $, $  {\mathtt y}^2 {\mathtt y}_{\mathtt x}  $ (and $ {\mathtt v}^3 $) 
there are no non-trivial periodic/quasi-periodic solutions, see Propositions \ref{nonex}-\ref{nonext}.
On the other hand, 
for  $  {\mathtt y}^3  $  the Birkhoff normal form is 
(partially) integrable by \cite{Po3} (for Dirichlet boundary conditions).  
\item {\bf Boundary conditions.} The  solutions of Theorem \ref{thm:DNLW} 
satisfy the Neumann boundary conditions $   {\mathtt y}_{\mathtt
x} (t,0) =  {\mathtt y}_{\mathtt x} (t, \pi )  = 0  $.
For proving the existence of solutions under Dirichlet boundary
conditions
it  would seem natural to substitute \eqref{REALITY} with 
 the oddness assumption  
\be\label{oddodd} 
g( - {\mathtt x}, -  {\mathtt y},  {\mathtt
y}_{\mathtt x},  {\mathtt v} ) = - g( {\mathtt x},   {\mathtt y},
{\mathtt y}_{\mathtt x}, {\mathtt v}  ) \, , 
\ee
so that 
the subspace of functions $ ({\mathtt y}, {\mathtt v} )({\mathtt x}) $ odd in $ {\mathtt x} $
is invariant under the flow evolution of \eqref{firstor}.  However, in order to find quasi-periodic solutions of \eqref{DNLW},
we need  the {\em real-coefficients property} \eqref{greco} which follows from \eqref{geven} and \eqref{REALITY},
 but not from \eqref{geven} and \eqref{oddodd}.   
It is easy to check that \eqref{geven}, \eqref{oddodd} and
\eqref{greco} imply the parity assumption \eqref{REALITY}. 
We decided to state the existence theorem in a form which requires the minimal assumptions.
Of course, if a nonlinearity satisfies  \eqref{geven}, \eqref{REALITY} and also 
\eqref{oddodd} we could look for quasi-periodic solutions satisfying Dirichlet boundary conditions. 
\item
{\bf Derivative vs quasi-linear NLW.}
It has been proved by
Klainermann-Majda \cite{KM}  that all classical solutions of Hamiltonian quasi-linear wave
equations like \be\label{QL} {\mathtt y}_{tt} = (1 + \sigma
({\mathtt y}_{\mathtt x})) {\mathtt y}_{{\mathtt x} {\mathtt x}}
\ee with $ \s^{(j)} (0) = 0 $, $ j = 1, \ldots, p -1 $, $
\sigma^{(p)} (0) \neq 0 $, 
do not admit smooth, small amplitude, periodic (a fortiori quasi-periodic)  solutions
except the constants.
Actually, any non constant solution of \eqref{QL},  with
sufficiently small initial data,
develops a singularity in finite time  in the second derivative $ {\mathtt y}_{  {\mathtt x} {\mathtt x}}  $. 
 In this respect \cite{KM} may suggest  that Theorem \ref{thm:DNLW} is 
 optimal regarding the order of (integer) derivatives in the nonlinearity.
Interestingly, the solutions of the derivative wave equation (which is a semilinear PDE) 
found in Theorem  \ref{thm:DNLW} are analytic in both time $ t $ and space $  \mathtt x $. 
Clearly the KAM approach developed in \cite{BBP3}  fails for
quasi-linear equations like \eqref{QL} because the auxiliary
vector field (whose flow defines the KAM transformations) is unbounded (of order $ 1$).
One could still ask for a KAM result for quasi-linear Klein Gordon equations (for which
Delort \cite{De1} proved  some steps of Birkhoff normal form). 
Note that   adding a mass term $ \mm \mathtt y $  in the left hand side of
\eqref{QL},  non constant 
periodic solutions of the form $\mathtt y(t, \mathtt x)=c(t)$ or ${\mathtt y}(t,\mathtt x)=c(\mathtt x)$ may occur.
\end{enumerate}

We finally complement the previous existence results with some negative results. 

\begin{proposition}\label{nonex} {\bf (\cite{BBP3})}
Let $ p \in \N $ be odd. The DNLW equations 
\be\label{NON1}
{\mathtt
y}_{tt} - {\mathtt y}_{{\mathtt x}{\mathtt x}} = {\mathtt
y}_{\mathtt x}^p + f({\mathtt y}) \, , \quad 
{\mathtt y}_{tt} - {\mathtt y}_{{\mathtt
x}{\mathtt x}} = \partial_{\mathtt x} ( {\mathtt y}^p ) +
f({\mathtt y}) \, , \quad {\mathtt x} \in {\T } \, , 
\ee
have no smooth quasi-periodic solutions except for trivial periodic
solutions of the form $\mathtt y(t,\mathtt x)=c(t)$. In particular
$f\equiv 0$ implies $c(t)\equiv const.$ 
\end{proposition}

This result is proved in \cite{BBP3} showing that 
$$
M({\mathtt y},  {\mathtt v} ) := \int_\T
{\mathtt y}_{\mathtt x} \, {\mathtt v} \, d{\mathtt x} 
$$
is a  Lyapunov function for \eqref{NON1}.
For wave equations, the
role of the  space variable $ {\mathtt x} $ and time variable $ t $ is
symmetric. A term like $ {\mathtt y}_t^p $ for an odd $ p $ is a
friction term which destroys the existence of quasi-periodic
solutions. Using 
$$
H({\mathtt y},  {\mathtt v} ) := \int_\T \frac{ {\mathtt v}^2}{2} +   \frac{{\mathtt y}_{\mathtt x}^2}{2}  -
F({\mathtt y}) \, d{\mathtt x}
$$
as a Lyapunov function we prove that: 

\begin{proposition}\label{nonext}{\bf (\cite{BBP3})}
Let $ p \in \N $ be odd. The DNLW equation
$$ 
{\mathtt y}_{tt} - {\mathtt y}_{{\mathtt x}{\mathtt x}} = {\mathtt y}_t^p +
f({\mathtt y}) \, , \quad {\mathtt x} \in \T  \, , 
$$ 
has no smooth quasi-periodic solutions except for trivial periodic
solutions of the form $\mathtt y(t,\mathtt x)=c(\mathtt x)$. In
particular $f\equiv 0$ implies $c(\mathtt x)\equiv const.$
\end{proposition}

The mass term $ \mm {\mathtt y} $ is, in general, 
necessary to have  existence of quasi-periodic solutions. The following non-existence results hold:

\begin{proposition}\label{FritzJ} {\bf (\cite{BBP3})}
The derivative NLW equation
\be\label{FJ}
{\mathtt y}_{tt} - {\mathtt y}_{{\mathtt x}{\mathtt
x}} 
= {\mathtt y}_t^2 \, , \quad {\mathtt x} \in \T \, , 
\ee
has no smooth solutions defined for all times except the constants. Moreover, for $ p , q \in \N $ even, 
\be\label{altre}
{\mathtt y}_{tt} - {\mathtt y}_{{\mathtt x}{\mathtt
x}} = {\mathtt y}_{\mathtt x}^{p} \, , 
\quad  
{\mathtt y}_{tt} - {\mathtt y}_{{\mathtt
x}{\mathtt x}} = {\mathtt y}_t^{p} \, , \quad
 {\mathtt y}_{tt} - {\mathtt y}_{{\mathtt
x}{\mathtt x}} = {\mathtt y}_{\mathtt x}^{p} + {\mathtt y}_t^{q}
\, , \quad {\mathtt x} \in \T \, , 
\ee
have no smooth periodic/quasi-periodic solutions except the constants.
\end{proposition}

\noindent
The blow-up result for \eqref{FJ} is proved by projecting the equation on the constants. 
The non-existence  results for \eqref{altre} may be obtained
simply by integrating the equations in $ (t, \mathtt x ) $.

\section{Ideas of proof: the abstract KAM theorem} 

The proof of Theorem \ref{thm:DNLW} is based on an abstract KAM
Theorem 
for reversible   infinite dimensional systems (Theorem 4.1 in \cite{BBP3})
which proves the existence of elliptic invariant tori and provides a reducible normal form 
around them. We now explain the main ideas and techniques of proof.

\smallskip

\noindent {\bf Complex formulation.} We extend \eqref{firstor} as
a first order system with complex valued variables $ ({\mathtt y},{\mathtt v}) \in
\C^n \times \C^n $. In  the unknowns
$$
u^+ :=\frac{1}{\sqrt 2} (D{\mathtt y}+\ii  {\mathtt v} )\,,\quad
u^- := \frac{1}{\sqrt 2} (D{\mathtt y}-\ii  {\mathtt v} )\,, \quad
D :=\sqrt{-\partial_{{\mathtt x}{\mathtt x}}+\mm} \, , \ \
\ii:=\sqrt{-1} \, ,
$$
system \eqref{firstor}  becomes the first order system
\be\label{brace}
\begin{cases}
u^+_t = - \ii D u^+ + \ii \mathtt g( u^+, u^- ) \cr u^-_t =  \ \
\ii D u^- -  \ii \mathtt g(u^+, u^-)
\end{cases}
\ee where
$$
\mathtt g(u^+, u^-) = \frac{1}{\sqrt{2}} \, g\Big( {\mathtt x},
D^{-1} \Big(\frac{ u^+ + u^- }{ \sqrt{2}}\Big) , D^{-1}
\Big(\frac{ u^+_{\mathtt x} +  u^-_{\mathtt x}}{  \sqrt{2}}\Big) ,
\frac{u^+ - u^-}{ \ii \sqrt{2}} \Big) \, .
$$
In \eqref{brace}, the dynamical variables  $ (u^+,  u^- ) $ are
independent. However, since $ g $ is real analytic (real on real),
the real subspace 
\be\label{subreal} 
{\mathtt R} := \{ \overline{ u^+} = u^- \} 
\ee
 is invariant under the flow evolution of \eqref{brace}, since
 \be\label{ronr}
\overline{\mathtt g(u^+ , u^- )} =\mathtt g (u^+ , u^- ) \, ,
\quad \forall (u^+ , u^- ) \in {\mathtt R} \, , \ee and the second
equation in \eqref{brace} reduces to the  complex conjugated of
the first one. Clearly, this corresponds to real valued solutions
$ ({\mathtt y}, {\mathtt v})  $ of the real system \eqref{firstor}. We say that system
\eqref{brace} is ``real-on-real".  
For systems satisfying this property it is
customary to use also the shorter notation
$$
(u^+, u^- ) = ( u, \bar u ) \, .
$$
Moreover the subspace of even functions 
\be\label{Esub} 
{\mathtt E} := \Big\{ u^+ ({\mathtt x}) = u^+ (-{\mathtt x}) \, , \ u^-
({\mathtt x}) = u^- (- {\mathtt x}) \Big\} 
\ee 
(see \eqref{evenx}) is invariant under the flow evolution of \eqref{brace}, by
\eqref{REALITY}. System \eqref{brace} is reversible with respect
to the involution 
\be\label{invco} 
S (u^+ , u^- ) := (u^-, u^+) \, , 
\ee 
(which is nothing but \eqref{invol1} in the variables $ (u^+,  u^- ) $).
\\[1mm] 
\noindent {\bf Dynamical systems formulation.} We
introduce infinitely many coordinates by Fourier transform
\be\label{Fourier} u^+ = \sum_{j \in \Z} u_j^+ e^{\ii j {\mathtt
x}} \, , \quad u^- = \sum_{j \in \Z}  u^-_j e^{- \ii j {\mathtt x}
}\,. \ee 
Then \eqref{brace} becomes the infinite dimensional
dynamical system
\begin{equation}\label{tuc}
\begin{cases}
{\dot u}^+_j = - \ii  \l_j u^+_j + \ii {\mathtt g}_j^+ ( \ldots, u^+_h, u^-_h, \ldots )   \cr 
{\dot u}^-_j = \ \ \ii  \l_j u^-_j - \ii {\mathtt g}_{j}^-  ( \ldots, u^+_h, u^-_h, \ldots )
\end{cases}  \quad \forall j \in \Z \, ,
\end{equation}
where \be\label{faramir}
 \l_j := \sqrt{j^2+\mm}
\ee are the eigenvalues of  $ D $ and
$$
\mathtt g_j^+ = \frac{1}{2 \pi} \int_{\T}\mathtt g
\Big( \sum_{h \in \Z} u_h^+ e^{\ii h {\mathtt x}} , \sum_{h \in
\Z} u^-_h e^{- \ii h {\mathtt x}} \Big) e^{- \ii j {\mathtt x}}
d{\mathtt x} \, ,  \quad \quad {\mathtt g}_j^-  := {\mathtt g}^+_{-j} \, .
$$
By \eqref{Fourier}, the ``real" subset $ {\mathtt R} $ in
\eqref{subreal} reads $ \overline{u_j^+} = u_j^- $
(this is the motivation for the choice of the signs in
\eqref{Fourier}) and, by \eqref{ronr},  the second equation
in \eqref{tuc} is the complex conjugated of the first one, namely
$$
\overline{ {\mathtt g}_j^+ } = {\mathtt g}_{j}^-  \quad {\rm when}
\quad  \overline{u_j^+} = u_j^- \, , \ \forall j \, .
$$
Moreover, the  invariant subspace $ \mathtt E $ of even functions in \eqref{Esub} reads, 
under Fourier transform, 
\be\label{parit} 
E := \Big\{ u^+_{j} =
u^+_{-j} \,  , \ u^-_{j} =  u^-_{-j} \, , \ \forall j \in \Z
\Big\} \, . 
\ee 
By \eqref{Fourier} the involution \eqref{invco} reads \be\label{Seven} S : (u^+_j,
u^-_j) \to (u^-_{-j}, u^+_{-j}) \, , \quad \forall j \in \Z \, .
\ee 
Finally, since $ g $ is real analytic, the assumptions
\eqref{geven} and \eqref{REALITY} imply the key property
\be\label{greco}
\mathtt g^\pm_j ( \ldots, u_j^+ , u^-_j, \ldots )  \ {\it has \ real \ Taylor \ coefficients} 
\ee 
in the variables $ (u_j^+ , u^-_j ) $.

\begin{remark}\label{rem:osc}
The previous property is compatible with oscillatory phenomena for
\eqref{tuc},  excluding friction phenomena. This is another strong
motivation for assuming \eqref{geven} and \eqref{REALITY}.
\end{remark}

\noindent {\bf Abstract KAM theorem.}
For every choice of   symmetric \textsl{tangential sites}
$$
 \mathcal I = \mathcal I^+\cup (-\mathcal I^+)
\quad {\rm with} \quad \mathcal I^+ \subset \N \setminus \{0\}  \,
, \ \sharp \mathcal I = n \, ,
$$ 
the linearized system \eqref{tuc} where $ \mathtt g_j^\pm  = 0 $ has the  invariant tori
$$
\Big\{ u_j\bar u_j=\xi_{|j|} > 0 ,\ {\rm for}\ j\in\mathcal{I}\,,\
u_j = \bar u_j = 0 \ {\rm for}\ j\not\in\mathcal{I}  \Big\}
$$
parametrized by the actions $ \xi = (\xi_j)_{j\in\mathcal{I}^+} $.
They correspond to the quasi-periodic solutions in
\eqref{ylineare}.

\smallskip

We first analyze the nonlinear dynamics of \eqref{tuc} close to the
origin, via a Birkhoff normal form reduction (see section 7 of \cite{BBP3}). 
This step depends on the nonlinearity $ g $  and
on the fact that the mass $ \mm >  0 $. Here we use \eqref{gleading}  to ensure that the third order  Birkhoff
normalized system is (partially) integrable and that the ``twist condition'' holds (the frequency-to-action map
is a diffeomorphism). 

Then we introduce action-angle coordinates on the tangential
variables: 
\be\label{AAV} u_j^+ = \sqrt{\xi_{|j|} + y_j} e^{\ii
x_j} \, , \, u^-_j = \sqrt{\xi_{|j|} + y_j} e^{-\ii x_j},  \  j
\in {\cal I} \, , \  (u_j^+, u^-_j)  = (z_j^+,z_j^-)\equiv(z_j, \bar z_j) \, , \, j
\notin {\cal I} \, , 
\ee where $ | y_j | < \xi_{|j|}  $.  Then,
system \eqref{tuc} transforms into  a parameter dependent family
of analytic systems of the form
\begin{equation}\label{HKAM}
\begin{cases}
{\dot x} = \oo( \xi ) +  P^{(x)} ( x, y, z, \bar z; \xi )\cr
 {\dot y} = \qquad \quad  P^{(y)} ( x, y, z, \bar z; \xi ) \cr {\dot z}_j
= - \ii \OO_j (\xi) z_j +  P^{(z_j)} (x,y,z,\bar z;\xi) \cr
{\dot {\bar z}}_j = \ \ \ii \O_j (\xi) {\bar z}_j + P^{({\bar
z}_j)}  (x,y,z,\bar z;\xi) \, , \ \ j \in \Z \setminus {\cal I} \, ,
\end{cases}
\end{equation}
where $ (x, y)  \in \T^n_s \times \C^n $,  $ z, \bar z  $ are
infinitely many variables, $ \oo(\xi)\in\R^n $,
$\OO(\xi)\in\R^\infty $. The frequencies $\oo_j(\xi),\, \OO_j (\xi)$
are close to the unperturbed frequencies $ \l_j $ in
\eqref{faramir} and satisfy $ \o_{-j} = \om_j $, $ \O_{-j} = \O_j $. 
The vector field $X$ in system \eqref{HKAM} has the form
$$
X = {\cal N} + P \, , 
$$
where  the normal form is 
\be\label{normalform}
{\cal N} := \oo( \xi ) \partial_x - \ii \OO_j (\xi) z_j \partial_{z_j} +  \ii \OO_j (\xi){\bar z}_j \partial_{\bar z_j} 
\ee
(we use the differential geometry notation for vector fields). The following properties hold: 
\begin{enumerate}
\item  {\sc reversible.} 
The vector field $ X = (X^{(x)}, X^{(y)}, X^{(z)}, X^{(\bar z)} ) $ 
is {\sc reversible}, namely 
$$
X \circ S = - S \circ X \, , 
$$
with respect to the involution 
$$ 
S : (x_j, y_j, z_j, \bar
z_j) \mapsto (-x_{-j}, y_{-j} , \bar z_{-j}, z_{-j}) \, , \
\forall j \in \Z \, , \quad S^2 = I \, , 
$$
 which is nothing but \eqref{Seven} in the variables \eqref{AAV}.
\item {\sc real-coefficients.} 
The components of the vector field
$$
X^{(x)} \, , \ \ii X^{(y)} \, , \  \ii X^{(z_j)}  \, , \ \ii X^{({\bar z}_j)}
$$
have {\it real} Taylor-Fourier
 coefficients  in the variables $ (x,y, z , \bar z ) $, see \eqref{greco}. 
\item {\sc real-on-real.} 
$$
X^{(x)}(v)= \overline{X^{(x)}(v)} \, , \
 X^{(y)}(v)=  \overline{X^{(y)}(v)} \, ,
\  X^{(z^-)}(v) =  \overline{X^{(z^+)}(v)} \, ,
\quad 
$$
for all $  v = (x,y,z^+ , z^-) $  such that $  x \in \T^n $, $ y \in \R^n $, $ \overline{z^+} = z^- $. 
\item {\sc even.}  The vector field $ X : E \to E $ where 
$$
 E := \Big\{ x_j = x_{-j}\,,\  y_j = y_{-j}\,, \ j \in \mathcal I\,,
 \  \ z_j = z_{-j}\,,\  \bar z_j = \bar z_{-j}\, , \ j \in \Z\setminus \mathcal I \Big\}
$$
is nothing but \eqref{parit} in the variables \eqref{AAV}. 
Hence the subspace $ E $ is invariant under the flow evolution of \eqref{HKAM}. 
\end{enumerate}

In  system \eqref{HKAM} we think $ x_j, y_j $, $ j \in {\cal I} $, $  z_{j}^\pm $, $ j \in \Z \setminus {\cal I} $, 
as {\it independent} variables and then we look for solutions
in the invariant subspace  $ E $, which means solutions of \eqref{brace} even in $ \mathtt x $. 
We proceed in this way because, in a phase space of functions even in $ \mathtt x $, 
the notion of {\it momentum}  (see   \eqref{momp}) 
is {\it not} well defined, as the Neumann boundary conditions break the translation invariance 
\footnote{In a more technical language 
we may see the above 
difficulty  by noting that, if 
$ z^\pm_{-j} \equiv z^\pm_j $,
the vector fields  $ z_{-j} \partial_{z_i} $ and $ z_j \partial_{z_i} $,  that 
have {\sc different} momentum, 
 would be identified.}. 
 On the other hand, the concept of  momentum 
is essential in order to verify the properties of  quasi-T\"oplitz vector fields, as explained after 
  \eqref{aweight}-\eqref{penalize}. 
That is why in \cite{BBP3} we  actually work in a phase space of $ 2 \pi $-periodic functions, 
where the notion   \eqref{momp}  of momentum is
well defined, and then we look for solutions in the invariant subspace \eqref{Esub} of even functions. 
\smallskip

A difficulty that arises  working in the whole phase space of periodic functions is that 
the linear frequencies $ \o_{-j} = \o_j  $, $ j \in {\cal I} $, $ \O_{-j} = \O_j $, $ j \in \Z \setminus {\cal I}$, 
 are {\it resonant}. Therefore, along the KAM iteration,
the monomial vector fields of the perturbation 
\begin{eqnarray*}
& & e^{\ii k\cdot x}\partial_{x_j} \, , \quad e^{\ii k\cdot x}y^i\partial_{y_j} \, ,  \  k \in \Z^n_{\rm odd} \, ,\,  \ \ |i| = 0,1, \  j \in \mathcal I \, ,  \\
& & e^{\ii k\cdot x} z_{\pm j} \partial_{z_j} \, , \quad  e^{\ii k\cdot x}  \bar z_{\pm j} \partial_{\bar z_j} , \ \ 
\forall  k \in \Z^n_{\rm odd}\,,\  j\in\Z\setminus \mathcal I\, ,
\end{eqnarray*}
where $ \Z^n_{\rm odd} := \{ k \in \Z^n : k_{-j} = - k_j \, , \forall j \in {\cal I} \} $, 
can {\it not} be averaged out. On the other hand, on the 
invariant subspace $ E $, where we look for the quasi-periodic solutions, the above terms can be replaced
by the constant coefficients monomial vector fields, obtained setting $ x_{-j} = x_j $, $ z^\pm_{-j} = z^\pm_j $.
More precisely,  we proceed as follows (section 5 of \cite{BBP3}):  we replace the nonlinear vector field
perturbation $ P $ with its symmetrized $ {\cal S} P $ defined, by linearity,  as
\begin{eqnarray}
&& \mathcal S(e^{\ii k\cdot x} \partial_{x_j}) :=
\partial_{x_j}\,, \quad  \mathcal S(e^{\ii k\cdot x}y^i\partial_{y_j}) :=
y^i\partial_{y_j}\,,\   \forall k \in \Z^n_{\rm odd} \, , \   |i| = 0,1, \  j \in \mathcal I \, ,  \nonumber \\
&& \mathcal S(e^{\ii k\cdot x} z_{\pm j} \partial_{z_j}) := z_{j} \partial_{z_j}
\,, 
\quad \mathcal S(e^{\ii k\cdot x}  \bar z_{\pm j} \partial_{\bar z_j}) :=  \bar z_{j} \partial_{\bar z_j}\,, \
\forall  k \in \Z^n_{\rm odd}\,,\ j \in\Z\setminus \mathcal I\, , \label{simmetriz}
\end{eqnarray}
and $ {\cal S} $ is the identity on the other monomial vector fields. 
Since
$$
P_{|E} = ({\cal S}P)_{|E }
$$
the two vector fields  determine the same dynamics on the invariant subspace $ E $ (Corollary 5.1 of \cite{BBP3}). 
Moreover $ {\cal S}P $ is reversible as well, and 
its weighted (see \eqref{aweight}) and quasi-T\"oplitz norms 
are (essentially) the same\footnote{This is due to the fact that the symmetrization procedure
in \eqref{simmetriz} does not increase the momentum, see  \eqref{momp}.} as those of $ P $  
(Proposition 5.2 of \cite{BBP3}). 

The homological equations  for a {\it symmetric} and {\it reversible} vector field perturbation 
can now be solved  (Lemma 5.1 of \cite{BBP3})
and  the remaining resonant term 
is a diagonal,
constant coefficients correction of the normal form \eqref{normalform} (also using the {\it real coefficients} property).
This procedure allows the KAM iteration to be carried out. Note
that, after this composite KAM step, the correction to the normal
frequencies described in \eqref{corrections} comes out from the symmetrized vector field  $ {\cal S} P $
and not  $ P $ itself. 

\smallskip

As in the Hamiltonian case \cite{BBP1}, a major difficulty of the
KAM iteration is to fulfill, at each iterative step, the second
order Melnikov non-resonance conditions. Actually, following the
formulation of the KAM theorem given in \cite{BBP1} it
is sufficient to verify
\begin{equation}\label{cervicale}
| \o^{\infty} (\xi) \cdot k +\O_i^{\infty}(\xi) -
\O_j^{\infty}(\xi)| \geq \frac{\g}{1+|k|^\tau} \, , \quad  \g > 0
\, ,
\end{equation}
only for the ``final" frequencies $ \om^\infty (\xi) $ and
$ \Om^\infty (\xi) $ and not along the inductive iteration.

As in \cite{BBP1} the key idea for verifying the second order
Melnikov non-resonance conditions \eqref{cervicale} for DNLW
is to  prove the higher order asymptotic decay estimate
\begin{equation}\label{delfi}
   \O_j^\infty(\xi) =
    j + a(\xi) +\frac{\mm}{2j}+  O (\frac{\g^{2/3}}{j})
    \quad {\rm for} \quad j \geq O(\g^{-1/3})
   \end{equation}
where $ a (\xi) $ is a constant independent of $ j $.

\smallskip

This property follows by introducing the notion of {\sl
quasi-T\"oplitz vector field}, see Definition 3.4 in \cite{BBP3}. The
new normal frequencies for a symmetric perturbation $ P = {\cal S}
P $ are $ \Om_j^+ = \Om_j + \ii P^{z_j,z_j}  $ where  the corrections $
P^{z_j,z_j} $ are the diagonal entries of the matrix defined by
\be\label{corrections}
P^{z,z}  z\partial_z  := \sum_{i,j} P^{z_i,z_j} z_j \partial_{z_i} \, , \quad  \ P^{z_i,z_j} := \int_{\T^n}
(\partial_{z_j}  P^{(z_i)})(x,0, 0, 0; \xi) \, d x \, .
\ee
Note that thanks to the {\it real-coefficients} property the corrections $ \ii P^{z_j,z_j} $ are real.
We say that a matrix $ P= P^{z,z}$ is quasi-T\"oplitz if it has the form
$$
P = T + R
$$
where $T$ is a T\"oplitz matrix (i.e. constant on the diagonals)
and $R$ is a ``small" remainder, satisfying in particular $ R_{jj} = O(1/ j) $.
Then  \eqref{delfi} follows with the constant $ a := T_{jj} $ which is
independent of $ j $.

The definition of quasi-T\"oplitz  vector field
is actually simpler than that of quasi-T\"oplitz function, used in the Hamiltonian context \cite{BBP1}, \cite{PX}. 
In turn, 
the notion of quasi-T\"oplitz function is weaker than the T\"oplitz-Lipschitz property, introduced 
by Eliasson-Kuksin in \cite{EK}. 

The quasi-T\"oplitz nature of the 
perturbation  is preserved along the  KAM iteration
(with slightly modified parameters) because 
the class of quasi-T\"oplitz vector fields is closed with respect to
\begin{enumerate}
\item Lie bracket $[ \ , \ ] $
(Proposition 3.1 of \cite{BBP3}),
\item Lie series
(Proposition 3.2 of \cite{BBP3}),
 \item
Solution of the homological equation
(Proposition 5.3 of \cite{BBP3}),
\end{enumerate}
namely the operations  which are used along the KAM iterative scheme.

\smallskip

An important difference with respect to  \cite{BBP1} 
is that we do {\it not} require the conservation of momentum, and so
our KAM theorem  applies 
to the DNLW equation \eqref{DNLW}  where the nonlinearity $ g $ may depend on the space
variable $ \mathtt x $. At a first insight this is a serious problem because
the properties of quasi-T\"oplitz functions as introduced in \cite{PX} and \cite{BBP1}, 
strongly rely on the conservation of  momentum.  

We remark that, anyway, the concept of {\it momentum} of a vector field is well defined (and useful) also 
if it is not conserved (prime integral), 
and so for PDEs 
which are not-invariant under $ \mathtt x $-translations (but with $ \mathtt x \in \T $). 
The momentum of a monomial vector field
$$
{\frak m}_{k,i,\a,\b; \mathtt v} := e^{\ii k \cdot x} y^i z^\a {\bar z}^\b \partial_{\mathtt v} \, ,  
\quad {\mathtt v} \in \{x, y, z_j, {\bar z}_j \} \, ,  
$$
(with  multi-indices notation
$ z^\a {\bar z}^\b := \Pi_{j \in \Z \setminus {\cal I}} \, z_j^{\a_j}  {\bar z}_j^{\b_j} $, $ \a_j $, $ \b_j \in \N $) 
is defined by
\begin{equation}\label{momp}
\pi(k,\a,\b;\vgot) :=
\left\{ \begin{array}{ll} \pi(k,\a,\b) &
\ \ \ \ {\rm if} \quad \vgot \in \{x_1,\ldots,x_n,y_1,\ldots,y_n\} \\
 \pi(k,\a,\b) -\s j & \ \  \ \ {\rm if} \quad \vgot = z_j^\s \, , \ \s = \pm \, , 
 \end{array}\right.
\end{equation}
where
$$
\pi(k,\a,\b):=
\sum_{i=1}^n \mathtt j_i k_i +\sum_{j \in \Z \setminus {\cal I}} (\a_j-\b_j)j \, , \quad 
$$
and $ { \cal I} := \{ {\mathtt j}_1, \ldots,  {\mathtt j}_n \}  $ are the tangential sites.
The monomial vector fields $ {\frak m}_{k,i,\a,\b; \mathtt v}  $ are nothing but 
the eigenvectors, with eigenvalues $ \ii \pi (k,\a,\b ; \mathtt v) $, 
of 
 the adjoint action of the momentum vector field
$ X_M := ( \mathtt j, 0, \ldots, \ii jz_j, \ldots , -\ii j \bar z_j, \ldots )$, namely
$$
[ {\frak m}_{k,i,\a,\b; \mathtt v} , X_M ] = \ii \pi (k,\a,\b ; \mathtt v){\frak m}_{k,i,\a,\b; \mathtt v} \, ,  
$$
see Lemma 2.1 of \cite{BBP3}.
 This is why it is convenient to use the 
 exponential basis in the Fourier decomposition \eqref{Fourier} instead of the cosine basis   $ \{ \cos (j \mathtt x) \}_{j \geq 0} $.

\begin{remark}
For a PDE which is translation invariant (namely the nonlinearity $ g $ is $ \mathtt x $-independent), 
all  the monomials  of the corresponding vector field $X$ have momentum equal to zero.
\end{remark} 

Then 
we overcome the impasse of the non-conservation of the momentum 
introducing the $ \mathtt a$-weighted majorant norm for
vector fields 
\be\label{aweight}
\norma X \norma_{s,r,\aaa} := 
\sup_{|y|_1 < r^2, \| z\|_{a,p}, \| \bar z\|_{a,p} < r}
\Big\| \Big(
 \sum_{k,i,\a, \b}  e^{\aaa|\pi(k,\a,\b;\vgot)|} | X_{k,i,\a, \b}^{(\vgot)}| e^{|k|s}
 |y^i| |z^\a| |{\bar z}^{\b}|
\Big)_{\vgot \in \mathtt{V}} \Big\|_{s,r} \, 
\ee
where $ r,s, \mathtt a >  0 $, $ {\mathtt V } := \{x,y, z_j, {\bar z}_j  \}, j \in \Z \setminus {\cal I} $, and 
$$
 \| (x,y,z, \bar z) \|_{s,r}   := \frac{|x|_\infty}{s} + \frac{|y|_1}{r^2} +  \frac{\| z\|_{a,p}}{r} +  \frac{\| \bar z\|_{a,p}}{r}  \, , 
 \ \ \| z \|_{a,p}^2 := \sum_{j \in \Z \setminus {\cal I}} |z_j|^2 e^{2 a |j|} \langle j \rangle^{2p} \, ,
$$ 
$ a \geq 0 $, $ p > 1 / 2 $ fixed (analytic spaces). The $ \| \ \|_{s,r,\aaa}$-norm  penalizes the
high-momentum monomials 
\be\label{penalize}
 \norma \Pi_{|\pi| \geq K} X \norma_{s,r,\aaa'}
 \leq  e^{-K(\aaa-\aaa')} \norma X\norma_{s,r,\aaa}  \,,\qquad 
 \forall\, 0\leq \aaa'\leq \aaa\, ,
\ee
so that only the {\it low-momentum} monomials vector
fields with $ |\pi| \leq K  $ are relevant (slightly decreasing $ \mathtt a' < \mathtt a $). 
This fact is crucial, in particular, in order to prove that the class of quasi-T\"oplitz  vector fields is closed
with respect to Lie brackets (Proposition 3.1 of \cite{BBP3}).



\noindent {\it Massimiliano Berti}, Dipartimento di Matematica e
Applicazioni ``R. Caccioppoli", Universit\`a degli Studi di Napoli
Federico II,  Via Cintia, Monte S. Angelo,
I-80126, Napoli, Italy,  {\tt m.berti@unina.it}. 
\\[0.5mm]
{\it Luca Biasco}, Dipartimento di Matematica, Universit\`a di
Roma 3, Largo San Leonardo Murialdo, I-00146, Roma, Italy, {\tt
biasco@mat.uniroma3.it}.
\\[0.5mm]
{\it Michela Procesi}, Dipartimento di Matematica ``G.
Castelnuovo", Universit\`a degli Studi di Roma la Sapienza,  P.le
A. Moro 5, 00185 ROMA, Italy,  {\tt
michela.procesi@mat.uniroma1.it}.
\\[0.5mm]
This research was supported by the European Research Council under
FP7 ``New Connections between dynamical systems and Hamiltonian
PDEs with small divisors phenomena" and partially by PRIN2009
grant ``Critical Point Theory and Perturbative Methods for
Nonlinear Differential Equations".

\end{document}